\documentclass[a4paper,11pt]{article}
\usepackage[latin1]{inputenc}
\usepackage[francais,english]{babel}    
\usepackage{amssymb}
\usepackage{amsmath}
\usepackage{amsthm}
\usepackage{multicol}
\usepackage{graphicx}
\usepackage{color}
\usepackage[T1]{fontenc}
\usepackage{yfonts}
\usepackage{placeins}
\usepackage{array}
\usepackage{dsfont}
\usepackage{stmaryrd}
\usepackage{verbatim}
\usepackage{caption}
\usepackage{transparent}

\newtheorem{thm}{Theorem}[section]
\newcommand{\bthm}{\begin{thm}}
\newcommand{\ethm}{\end{thm}}

\newtheorem{thmi}{Theorem}
\newcommand{\bthmi}{\begin{thmi}}
\newcommand{\ethmi}{\end{thmi}}

\newtheorem{cori}[thmi]{Corollary}
\newcommand{\bcori}{\begin{cori}}
\newcommand{\ecori}{\end{cori}}

\newtheorem{mthm}{Theorem}
\newcommand{\bmthm}{\begin{mthm}}
\newcommand{\emthm}{\end{mthm}}
\newtheorem{mcor}[mthm]{Corollary}
\newcommand{\bmcor}{\begin{mcor}}
\newcommand{\emcor}{\end{mcor}}

\newtheorem*{conj}{Conjecture}
\newcommand{\bconj}{\begin{conj}}
\newcommand{\econj}{\end{conj}}

\newtheorem*{question}{Question}
\newcommand{\bq}{\begin{question}}
\newcommand{\eq}{\end{question}}

\newtheorem*{thn}{Theorem}
\newcommand{\bthn}{\begin{thn}}
\newcommand{\ethn}{\end{thn}}

\newtheorem{exo}{Exercise}
\newcommand{\bex}{\begin{exo}}
\newcommand{\eex}{\end{exo}}

\newtheorem{sol}{Solution}
\newcommand{\bsol}{\begin{sol}}
\newcommand{\esol}{\end{sol}}

\newtheorem{pro}[thm]{Proposition}
\newcommand{\bpro}{\begin{pro}}
\newcommand{\epro}{\end{pro}}

\newtheorem{cor}[thm]{Corollary}
\newcommand{\bcor}{\begin{cor}}
\newcommand{\ecor}{\end{cor}}

\newtheorem{lem}[thm]{Lemma}
\newcommand{\blem}{\begin{lem}}
\newcommand{\elem}{\end{lem}}

\theoremstyle{definition}

\newtheorem{defi}[thm]{Definition}
\newcommand{\bdf}{\begin{defi}}
\newcommand{\edf}{\end{defi}}

\newtheorem*{defis}{Definition}
\newcommand{\bdfs}{\begin{defis}}
\newcommand{\edfs}{\end{defis}}

\newtheorem*{rmk}{Remark}
\newcommand{\brk}{\begin{rmk} \upshape}
\newcommand{\erk}{\end{rmk}}

\newtheorem*{exe}{Example}
\newcommand{\bexe}{\begin{exe} \upshape}
\newcommand{\eexe}{\end{exe}}

\newtheorem*{pre}{Proof}
\newcommand{\bp}{\begin{pre} \upshape}
\newcommand{\ep}{\hfill \qed \end{pre}}
\newcommand{\epp}{\end{pre}}

\newtheorem{mainthm}{Theorem}

\newtheorem{mainpro}[mainthm]{Proposition}

\newcommand{\beq}{\begin{eqnarray*}}
\newcommand{\eeq}{\end{eqnarray*}}

\newcommand{\beqn}{\begin{equation}}
\newcommand{\eeqn}{\end{equation}}

\newcommand{\ben}{\begin{enumerate}}
\newcommand{\een}{\end{enumerate}}
\newcommand{\bit}{\begin{itemize} \renewcommand{\labelitemi}{$\bullet$} \renewcommand{\labelitemii}{$\star$}}
\newcommand{\eit}{\end{itemize}}

\newcommand{\bfg}{
\begin{figure}[H]
\begin{center}}
\newcommand{\efg}{
\end{center}
\end{figure}
\FloatBarrier}

\newcolumntype{M}[1]{>{\raggedright}m{#1}}

\newcommand{\R}{\mathbb{R}}

\newcommand{\N}{\mathbb{N}}
\newcommand{\Z}{\mathbb{Z}}
\newcommand{\C}{\mathbb{C}}
\newcommand{\B}{\mathbb{B}}

\newcommand{\K}{\mathbb{K}}

\renewcommand{\SS}{\mathbb{S}}

\newcommand{\bs}{\symbol{92}}

\newcommand{\ov}{\overline}

\newcommand{\llb}{\llbracket}
\newcommand{\rrb}{\rrbracket}

\newcommand{\st}{\, | \,}

\newcommand{\ra}{\rightarrow}

\newcommand{\liml}{\lim\limits}

\renewcommand{\geq}{\geqslant}
\renewcommand{\leq}{\leqslant}

\newcommand{\SL}{\operatorname{SL}}
\newcommand{\SO}{\operatorname{SO}}
\newcommand{\SU}{\operatorname{SU}}

\newcommand{\pif}{{+\infty}}

\newcommand{\sign}{\begin{flushright}
Thomas Haettel \\
Universit\'e de Montpellier \\
Institut Montpelli\'erain Alexander Grothendieck \\
CC051 \\
Place Eug\`ene Bataillon \\
34095 Montpellier Cedex 5 \\
France \\
thomas.haettel@math.univ-montp2.fr
\end{flushright}}

\makeatletter
\def\Ddots{\mathinner{\mkern1mu\raise\p@
\vbox{\kern7\p@\hbox{.}}\mkern2mu
\raise4\p@\hbox{.}\mkern2mu\raise7\p@\hbox{.}\mkern1mu}}
\makeatother

\makeatletter
\def\maketitles{%
  \null
  \thispagestyle{empty}%
  \vfill
  \begin{center}\leavevmode
    \normalfont
    {\LARGE \@title\par}%
    \vskip 1.2cm
    {\large \@author\par}%
    \vskip 1.2cm
    {\large \@subtitle\par}%
    \vskip 0.8cm
    {\large \@date\par}%
  \end{center}%
  \vfill
  \null
  \cleardoublepage
  }

\def\date#1{\def\@date{#1}}
\def\author#1{\def\@author{#1}}
\def\title#1{\def\@title{#1}}
\def\subtitle#1{\def\@subtitle{#1}}
\makeatother

\newtheorem*{procube}{Proposition~\ref{pro:embedding_cube}}

\title{Higher rank lattices are not coarse median}
\author{Thomas Haettel}
\date{\today}

\begin{document}

\selectlanguage{english}

\maketitle

\begin{center}
\begin{minipage}{0.8\textwidth}
\textsc{Abstract.} We show that symmetric spaces and thick affine buildings which are not of spherical type $A_1^r$ have no coarse median in the sense of Bowditch. As a consequence, they are not quasi-isometric to a CAT(0) cube complex, answering a question of Haglund. Another consequence is that any lattice in a simple higher rank group over a local field is not coarse median.
\end{minipage}
\end{center}

\let\thefootnote\relax

\section*{Introduction}

A metric space $(X,d)$ is called metric median if for each $(x,y,z) \in X^3$, the three intervals $I(x,y)$, $I(y,z)$ and $I(x,z)$ intersect in a single point, where $I(x,y)=\{p \in X \st d(x,p)+d(p,y)=d(x,y)\}$. This point is called the median of $x$, $y$ and $z$. The rank of $(X,d)$ is then defined as the maximal dimension $r$ of an embedded cube $\{0,1\}^r$. The relationship between groups and median metric spaces is rich and has been studied through many points of view: Haagerup property, property (T), action on a CAT(0) cube complex, action on a space with (measured) walls... (see~\cite{chepoi_median_cat0}, \cite{chatterji_drutu_haglund}, \cite{chatterji_fernos_iozzi}, \cite{chatterji_niblo}, \cite{bowditch_embedding_trees}, \cite{bowditch_median_metric}...).

Bowditch recently introduced the notion of a coarse median on a metric space (see~\cite{bowditch_coarse_median}), in order to gather in the same setting CAT(0) cube complexes and Gromov hyperbolic spaces. A metric space is Gromov-hyperbolic if and only if every finite subset admits a good metric comparison with a tree (see for instance~\cite[Theorem~12, p.~33]{ghys_delaharpe}). Bowditch's definition of a coarse median is having a good metric comparison of every finite subset with a metric median space, or equivalently with a CAT(0) cube complex according to Chepoi (see~\cite{chepoi_median_cat0}).

\bdfs[Bowditch] Let $(X,d)$ be a metric space. A map $\mu:X^3 \ra X$ is called a \emph{coarse median} if there exist $k \in [0,\pif)$ and $h:\N \ra [0,\pif)$ such that
\bit
\item For all $a,b,c,a',b',c' \in X$, we have$$d(\mu(a,b,c),\mu(a',b',c')) \leq k(d(a,a')+d(b,b')+d(c,c')) + h(0).$$
\item For each finite non-empty set $A \subset X$, with $|A| \leq p$, there exists a finite median algebra $(\Pi,\mu_\Pi)$ and maps $\pi:A \ra \Pi$, $\lambda:\Pi \ra X$ such that for every $x,y,z \in X$, we have
$$d(\lambda \mu_\Pi(x,y,z),\mu(\lambda x,\lambda y,\lambda z)) \leq h(p),$$
and for every $a \in A$, we have $d(a,\lambda(\pi(a))) \leq h(p)$.
\eit
If $\Pi$ can be chosen (independently of $p$) to have rank at most $r$, we say that $\mu$ has rank at most $r$.
\edfs

A finitely generated group is said to be coarse median if some Cayley graph has a coarse median (not necessarily equivariant under the group action). Bowditch showed that a coarse median group is finitely presented, and has at most quadratic Dehn function (see~\cite[Corollary~8.3]{bowditch_coarse_median}). Furthermore, Chatterji-Ruane's criterion (see~\cite{chatterji_ruane}) applies to show that a coarse median group of finite rank has property (RD) (see~\cite[Theorem~9.1]{bowditch_embedding_trees}). Moreover, if a group has a coarse median of rank at most $r$, there is no quasi-isometric embedding of $\R^{r+1}$ into that group. Bowditch also showed that the existence of a coarse median is a quasi-isometry invariant, that a group is Gromov hyperbolic if and only if it is coarse median of rank $1$, and that a group hyperbolic relative to coarse median groups is itself coarse median (see~\cite{bowditch_rh_coarse_median}). Furthermore, Bowditch showed that the mapping class group of a surface of genus $g$ with $p$ punctures is coarse median of rank $3g-3+p$, hereby recovering Behrstock and Minsky's result that the mapping class group has property (RD) (see~\cite{behrstock_minsky_centroids}), and the rank theorem (see~\cite{hamenstadt_qirigidity} and \cite{behrstock_minsky_rank}).

Since most known examples of coarse median groups have some non-positive curvature features, Bowditch asked in~\cite{bowditch_coarse_median} whether higher rank symmetric spaces, or even CAT(0) spaces, admit coarse medians. In this article, we give a negative answer to this question.

\begin{mainthm} \label{thm:no median on Xinfty} Let $X$ be a thick affine building of spherical type different from $A_1^r$. There is no locally convex Lipschitz median on $X$. \end{mainthm}

By considering asymptotic cones and using work of Kleiner-Leeb 	and Bowditch, we deduce the following.

\begin{mainthm} \label{thm:no_cm_on_X_maindirection} Let $X$ be a symmetric space of non-compact type, or a thick affine building, of spherical type different from $A_1^r$. Then $X$ has no coarse median. \end{mainthm}

A consequence of this result is the classification of symmetric spaces of non-compact type and affine buildings which are coarse median.

\begin{mainthm} \label{thm:no_cm_on_X} Let $X$ be a symmetric space of non-compact type, or a thick affine building. There exists a coarse median on $X$ if and only if the spherical type of $X$ is $A_1^r$. \end{mainthm}

Note that the coarse median is not assumed to be equivariant by any group.

\medskip

Haglund asked if a higher rank symmetric space or affine building is quasi-isometric to a CAT(0) cube complex, and we give a negative answer:

\begin{mainthm} \label{thm:not_qi_CCC} Let $X$ be a symmetric space of non-compact type, or a discrete, thick affine building. Then $X$ is quasi-isometric to a CAT(0) cube complex if and only if the spherical type of $X$ is $A_1^r$ \end{mainthm}

Note that the CAT(0) cube complex is not assumed to be of finite dimension, and it could also be endowed with the $L^p$ distance for any $p \in [1,\infty]$.

Also note that Theorem~\ref{thm:not_qi_CCC} still holds if we consider non-discrete thick affine buildings and non-discrete CAT(0) cube complexes.

\medskip

Furthermore, for uniform lattices in semisimple Lie groups, Property (RD) implies the Baum-Connes conjecture without coefficient (see~\cite{lafforgue_RD}). Property (RD) has been proved notably for uniform lattices in $\SL(3,\K)$, where $\K$ is a local field (see~\cite{ramagge_robertson_steger}, \cite{lafforgue_SL3} and \cite{chatterji_SL3}). Valette conjectured that uniform lattices in semisimple Lie groups satisfy property (RD). Since being coarse median implies property (RD), one could ask if this could be a way to prove property (RD) for higher rank uniform lattices. Even though it might follow from~\cite{ramagge_robertson_steger} that looking only at coarse medians is not enough for $\SL(3,\K)$, the following makes it clear.

\begin{mainthm} \label{thm:no_cm_on_G} Let $\K$ be a local field, let $G$ be the group of $\K$-points of a simple algebraic group without compact factors and let $\Gamma$ be a lattice in $G$. If $\Gamma$ is coarse median, then $G$ has $\K$-rank $1$.
\end{mainthm}

Note that, due to property (T), higher rank lattices do not admit unbounded actions on median metric spaces (see~\cite{chatterji_drutu_haglund}). But in the coarse median setting this is not a consequence of property (T), since for instance every hyperbolic group with property (T) is coarse median.

\medskip

In the $\K$-rank $1$ case, finding which non-uniform lattices are coarse median is harder. Here we summarize what is known.

\begin{mainpro} \label{pro:rank1} Let $\K$ be a local field, let $G$ be the group of $\K$-points of a simple algebraic group without compact factors of $\K$-rank $1$, and let $\Gamma$ be a lattice in $G$.
\bit
\item If $\Gamma$ is uniform in $G$, then $\Gamma$ is coarse median.
\item If $G$ is locally isomorphic to $\SO_0(n,1)$ for some $n \geq 2$, then $\Gamma$ is coarse median.
\item If $G$ is locally isomorphic to $\SU(2,1)$, then $\Gamma$ is not coarse median.
\eit
\end{mainpro}

\medskip

In the proof of Theorem~\ref{thm:no median on Xinfty}, we establish the following result, which is of independent interest.

\begin{mainpro} \label{pro:embedding_cube} Let $X$ be a connected metric space, with a Lipschitz locally convex median of rank $r$. There exists a median, biLipschitz embedding of the $r$-cube $[0,1]^r$ into $X$ with convex image. \end{mainpro}

\textbf{Organization of the paper:} In Section~\ref{sec:definitions}, we recall the general definitions of median algebras. In Section~\ref{sec:ultralimit}, we recall Kleiner-Leeb's and Bowditch's work on asymptotic cones, and we prove that Theorem~\ref{thm:no median on Xinfty} implies Theorems~\ref{thm:no_cm_on_X_maindirection} and \ref{thm:no_cm_on_X}. 

Section~\ref{sec:convex_apt} and \ref{sec:gate} are devoted to the proof of Theorem~\ref{thm:no median on Xinfty}. We consider a thick affine building $X$ which has a locally convex Lipschitz median. In Section~\ref{sec:convex_apt} we prove Proposition~\ref{pro:embedding_cube}, which provides us with a convex cube in $X$. In Section~\ref{sec:gate}, by considering a tangent cone of $X$ in the cube we can assume that some apartment $F$ of $X$ is isomorphic to a vector space with the standard $L^1$ median. Considering singular geodesics in $F$, we prove that $X$ has spherical type $A_1^r$.

Finally in Section~\ref{sec:corollaries}, we prove the main consequences of Theorem~\ref{thm:no_cm_on_X}, which are Theorem~\ref{thm:not_qi_CCC}, Theorem~\ref{thm:no_cm_on_G} and Proposition~\ref{pro:rank1}.

\medskip

\textbf{Acknowledgments:} The author would like to thank very warmly Brian Bowditch, who motivated this work and answered many questions in a precise way. The author would like to thank the anonymous referee, who made insightful comments and explained the subtility of the rank one case. The author would also like to thank Linus Kramer, Frédéric Haglund and Yves de Cornulier, who answered questions about this work. The author would also like to thank Samuel Tapie and Jean Lécureux, for interesting discussions at an early stage of this work.

\section{Median algebras}

\label{sec:definitions}

\bdf Let $X$ be a set. A map $\mu:X^3 \ra X$ is called a \emph{median} on $X$ if for all $a,b,c,d,e$ in $X$, it satisfies
\begin{eqnarray} & \star &\mu(a,b,c)=\mu(b,a,c)=\mu(b,c,a) \mbox{ ($\mu$ is symmetric)}, \nonumber \\
& \star &\mu(a,a,b)=a, \nonumber \\
\label{eq:abcde} & \star &\mu(a,b,\mu(c,d,e))=\mu(\mu(a,b,c),\mu(a,b,d),e). \end{eqnarray}
The pair $(X,\mu)$ is called a \emph{median algebra}.
\edf

\brk\
\bit
\item There is a unique median on the set $\{0,1\}$.
\item We can consider the product median on the $n$-cube $\{0,1\}^n$.
\eit
\erk

\bdf Let $(X,\mu)$, $(X',\mu')$ be median algebras. A map $f:X \ra X'$ is called \emph{median} if for every $x,y,z \in X$, we have $\mu'(f(x),f(y),f(z))=f(\mu(x,y,z))$. If furthermore $f$ is injective, it is called a \emph{median embedding}. \edf

\bdf Let $(X,\mu)$ be a median algebra. If every median embedding of an $n$-cube $\{0,1\}^n \ra X$ satisfies $n \leq r$, we say that $X$ has \emph{rank} at most $r$. \edf

\bdf Let $(X,\mu)$ be a median algebra. If $a,b \in X$, the \emph{interval} between $a$ and $b$ is $I(a,b)=\{c \in X \st \mu(a,b,c)=c\}$. A subset $C \subset X$ is called \emph{convex} if for every $a,b \in C$, we have $I(a,b) \subset C$. \edf

If $(X,d)$ is a metric space, a weakening of the notion of metric median is the following.

\bdf Let $(X,d)$ be a metric space. An abstract median $\mu$ on $X$ is called
\bit
\item continuous if $\mu:X^3 \ra X$ is a continuous map,
\item \emph{Lipschitz} if there exists a constant $k \in [0,\pif)$ such that $\mu$ is $k$-Lispchitz with respect to each variable, i.e. for every $a,b,c,a',b',c' \in X$, we have
$$d(\mu(a,b,c),\mu(a',b',c')) \leq k(d(a,a')+d(b,b')+d(c,c')),$$
\item \emph{locally convex} if each point of $X$ has a basis of neighbourhoods consisting of convex subsets of $X$.
\eit
 \edf

Here is an example of a continuous median on $\R^2$ which is not Lipschitz: consider the image $\mu$ of the standard $L^1$ median by some non-Lipschitz diffeomorphism of $\R^2$. If we consider $\R^2$ with the standard $L^1$ distance and the new median $\mu$, then $\mu$ is a continuous (even differentiable) median, but it is not Lipschitz.

\bdf Let $(X,d)$ be a metric space, let $\mu$ be a continuous median on $X$, and let $C \subset X$ be a non-empty closed, locally compact convex subset of $X$. Then for each $x \in X$, there exists a unique $\pi_C(x) \in C$, called the \emph{gate projection} of $x$ onto $C$, such that for every $y \in C$, we have $\pi_C(x) \in I(x,y)$. The map $\pi_C : X \ra C$ is called the gate projection, it is a continuous map. If $\mu$ is $k$-Lipschitz, then $\pi_C$ is a $k$-Lipschitz map. \edf

Now we recall the definition of walls in a median algebra.

\bdf Let $(X,\mu)$ be a median algebra. A \emph{wall} in $X$ is a pair $W=\{H^+(W),H^-(W)\}$, where $H^+(W)$ and $H^-(W)$ are non-empty convex disjoints subsets of $X$ whose union is equal to $X$. \edf

\blem \cite[Lemma~6.1]{bowditch_coarse_median} \label{lem:exists wall} Let $(X,\mu)$ be a median algebra, and let $A,B$ be disjoint convex subsets of $M$. There exists a wall $W=\{H^+(W),H^-(W)\}$ in $X$ separating $A$ from $B$, i.e. such that $A \subset H^{\pm}(W)$ and $B \subset H^{\mp}(W)$. \elem

\blem \cite[Lemma~7.3]{bowditch_coarse_median} \label{lem:exists strong wall} Let $(X,d)$ be a metric space, and let $\mu$ be a continuous locally convex median on $X$. Let $a,b$ be distinct points of $X$. There exists a wall $W=\{H^+(W),H^-(W)\}$ in $X$ strongly separating $a$ from $b$, i.e. such that $a \in X \bs \ov{H^-(W)}$ and $b \in X \bs \ov{H^+(W)}$. \elem

\blem \label{lem:rank_lw} \cite[Lemma~7.5]{bowditch_coarse_median} Let $X$ be a metric space, and let $\mu$ be a continuous locally convex median on $X$. For each wall $W=\{H^+(W),H^-(W)\}$ in $X$, the subset $L(W) = \ov{H^+(W)} \cap \ov{H^-(W)}$ is a convex median subalgebra of $X$, of rank at most $r-1$ if the rank of $\mu$ is $r$. \elem

\section{Ultralimits of spaces and coarse medians}

\label{sec:ultralimit}

In~\cite{kleiner_leeb}, Kleiner and Leeb developed a geometric definition of spherical and affine buildings, and in particular they studied their asymptotic cones.

\bthm \cite[Theorem~1.2.1]{kleiner_leeb} \label{thm:kleiner_leeb_asymptotic_cone} Let $X$ be a symmetric space of non-compact type or a thick affine building. Then any asymptotic cone of $X$ is a thick affine building of the same spherical type as $X$. \ethm

They also proved that any tangent cone of an affine building is an affine building:

\bthm \cite[Theorem~5.1.1]{kleiner_leeb} \label{thm:ultralimit of affine building} Let $(X,d)$ be an affine building, let $\omega$ be a non-principal ultrafilter on $\N$, let $(x_n)_{n \in \N}$ be a sequence in $X$ and let $(\lambda_n)_{n \in \N}$ be a sequence in $(0,\pif)$ such that $\liml_{\omega} \lambda_n = \pif$. Let $(X_\infty,d_\infty,x_\infty)$ be the $\omega$-ultralimit of $(X_n,\lambda_n d,x_n)$. Then $(X_\infty,d_\infty)$ is an affine building. Furthermore, if $X$ is thick, then $X_\infty$ is thick. The affine Weyl group of $X_\infty$ acts transitively on each apartment of $X_\infty$. \ethm

One motivation for Bowditch's definition of coarse median is that it behaves well when one considers asymptotic cones.

\bthm \cite[Theorem~2.3]{bowditch_coarse_median} \label{thm:limit median} Let $(X,d)$ be a metric space, and let $\mu$ be a $(k,h)$-coarse median on $X$. Then on any asymptotic cone $(X_\infty,d_\infty)$ of $(X,d)$, $\mu$ defines a $k$-Lipschitz, locally convex median $\mu_\infty$. \ethm

We can now prove that Theorem~\ref{thm:no median on Xinfty} implies Theorem~\ref{thm:no_cm_on_X_maindirection}.

\bp Let $X$ be a symmetric space of non-compact type, or a thick affine building, of spherical type different from $A_1^r$. By Theorem~\ref{thm:kleiner_leeb_asymptotic_cone} any asymptotic cone $X_\infty$ is a thick affine building. If there existed a coarse median $\mu$ on $X$, it would give rise by Theorem~\ref{thm:limit median} to a locally convex Lipschitz median on $X_\infty$, which contradicts Theorem~\ref{thm:no median on Xinfty} since the spherical type of $X_\infty$ is not $A_1^r$. Hence there is no coarse median on $X$.\ep 

We can also prove that Theorem~\ref{thm:no_cm_on_X_maindirection} implies Theorem~\ref{thm:no_cm_on_X}.

\bp Let $X$ be a symmetric space of non-compact type or an affine building of spherical type $A_1^r$. If $X$ is a symmetric space, it is a product of rank $1$ symmetric spaces, which are Gromov hyperbolic, so $X$ has a coarse median. If $X$ is an affine building, if we endow it with the $L^1$ metric it becomes a metric median space. In particular, this median is a coarse median with respect to any usual metric on $X$, which is equivalent to the $L^1$ metric.\ep 

\section{Existence of a convex cube}

\label{sec:convex_apt}

In this Section, we will prove Proposition~\ref{pro:embedding_cube}, which we recall here.

\begin{procube} Let $X$ be a connected metric space, with a Lipschitz locally convex median of rank $r$. There exists a median, biLipschitz embedding of the $r$-cube $[0,1]^r$ into $X$ with convex image. \end{procube}

Fix $X$ a geodesic metric space, and $\mu : X^3\ra X$ a Lipschitz median.

\bdf A continuous path $p:I \ra X$, where $I\subset \R$ is an interval, is said to be \emph{monotone} is for each $t_1 < t_2 < t_3$ in $I$, we have $\mu(p(t_1),p(t_2),p(t_3))=p(t_2)$. \edf
%
%

To prove Proposition~\ref{pro:embedding_cube}, we will need the following two Lemmas.

\blem \label{lem:gate cube} Let $X$ be a connected metric space, with a continuous locally convex median $\mu$, and let $f:\{0,1\}^r \ra X$ be a median embedding of the $r$-cube, and let $W$ be a wall in $X$ strongly separating $f(0,\ldots,0)$ and $f(1,0,\ldots,0)$. There exists a median embedding $g:\{0,1\}^r \ra X$ such that for every $t \in \{0\} \times \{0,1\}^{r-1}$, we have $g(t)=f(t)$ and for every $t \in \{1\} \times \{0,1\}^{r-1}$, we have $g(t) \in L(W)$. \elem

\bp Note that if we knew that $L(W)$ was locally compact, projecting the half-cube $\{1\} \times \{0,1\}^{r-1}$ using the gate projection onto $L(W)$ would immediately give the result.

Since intervals are connected, we can consider $a \in I(f(0,\ldots,0),f(1,0,\ldots,0)) \cap L(W)$. Define
\beq g : \{0,1\}^r & \ra & X \\
t \in \{0\} \times \{0,1\}^{r-1} & \mapsto & f(t) \\
t \in \{1\} \times \{0,1\}^{r-1} & \mapsto & \mu(f(0,t_2,\ldots,t_r),a,f(t)). \eeq
Since $L(W)$ is convex, we deduce that for evert $t \in \{1\} \times \{0,1\}^{r-1}$, we have $g(t) \in L(W)$. 

Using repeatedly Property~(\ref{eq:abcde}), we prove that $g$ is a median map. As a consequence, if for some $t,t' \in \{0,1\}^r$ we have $g(t)=g(t')$, then
\beq f(0,t_2,\ldots,t_r) &=& \mu(g(0,t_2,\ldots,t_r),g(t),g(0,t'_2,\ldots,t'_r)) \\
&=& \mu(g(0,t_2,\ldots,t_r),g(t'),g(0,t'_2,\ldots,t'_r)) \\
&=& f(0,t'_2,\ldots,t'_r),\eeq
hence $(0,t_2,\ldots,t_r)=(0,t'_2,\ldots,t'_r)$, so $t=t'$: $g$ is injective.
\ep

\blem \label{lem:convex cube} Let $(X,\mu)$ be a median algebra. Assume there exists a median embedding of the $r$-cube $f: [0,1]^r \ra X$, such that the image by $f$ of any edge of $[0,1]^r$ is convex. Then the image of $f$ is convex in $X$. \elem

\bp
For each $k \in \llb 1,r \rrb$, let $e_k=(0,\ldots,0,1,0,\ldots,0)$ (the $1$ is in the $k^{th}$ position). Let $x \in I(f(0),f(e_1+\ldots+e_r))$. For each $k \in \llb 1,r \rrb$, since the image by $f$ of the edge $[0,e_k]$ is convex, we deduce that $I(f(0),f(e_k))=f([0,e_k])$, so there exists $t_k \in [0,1]$ such that $\mu(f(0,\ldots,0),x,f(e_k))=f(t_ke_k)$. We will show by induction on $k \in \llb 0,r\rrb$ that $f(t_1e_1+\ldots+t_ke_k)=\mu(f(0),x,f(e_1+\ldots+e_k))$.

For $k=0$ this is immediate using Property~(\ref{eq:abcde}), so assume that for some $k<r$ we have $f(t_1e_1+\ldots+t_ke_k)=\mu(f(0),x,f(e_1+\ldots+e_k))$. Then
\beq f(t_1e_1+\ldots+t_{k+1}e_{k+1}) &=& \mu(f(t_1e_1+\ldots+t_ke_k),f(t_{k+1}e_{k+1}),f(e_1+\ldots+e_{k+1})) \\
&=& \mu(\mu(f(0),x,f(e_1+\ldots+e_k)),\mu(f(0),x,f(e_{k+1})),f(e_1+\ldots+e_{k+1})) \\
&=& \mu(f(0),x,\mu(f(e_1+\ldots+e_k),f(e_{k+1}),f(e_1+\ldots+e_{k+1})) \\
&=& \mu(f(0),x,f(e_1+\ldots+e_{k+1})).\eeq

As a consequence, for $k=r$ we deduce that
$$f(t_1e_1+\ldots+t_re_r)=\mu(f(0),x,f(e_1+\ldots+e_r))=x,$$
as $x \in I(f(0),f(e_1+\ldots+e_r))$. So we have proved that the image of $f$ is equal to the interval $I(f(0),f(e_1+\ldots+e_r))$, which is convex.
\ep

We can now prove Proposition~\ref{pro:embedding_cube}.

\bp[Proof of Proposition~\ref{pro:embedding_cube}]
Since the rank of the median $\mu$ is $r$, consider a median embedding $f:\{0,1\}^r \ra X$. Applying $2r$ times Lemma~\ref{lem:gate cube}, up to replacing $f$ by another median embedding of $\{0,1\}^r$ into $X$, we can assume that for each $i \in \llb 1,r\rrb$ and $\varepsilon \in \{0,1\}$, the image under $f$ of the codimension $1$ face $\{0,1\}^{i-1} \times \{\varepsilon\} \times \{0,1\}^{r-1-i}$ is included in a closed convex subspace $L(W_{i,\varepsilon})$ of $X$, where $W_{i,\varepsilon}$ is a wall of $X$. According to Lemma~\ref{lem:rank_lw}, each $L(W_{i,\varepsilon})$ has rank at most $r-1$, and since it contains the image by $f$ of the $(r-1)$-cube $\{0,1\}^{i-1} \times \{\varepsilon\} \times \{0,1\}^{r-1-i}$, we deduce that each $L(W_{i,\varepsilon})$ has rank $r-1$.

For $i,j \in \llb 1,r\rrb$ distinct and $\varepsilon,\varepsilon' \in \{0,1\}$, since $L(W_{i,\varepsilon}) \cap L(W_{j,\varepsilon'}) = L(L(W_{i,\varepsilon}) \cap W_{j,\varepsilon'})$, where $L(W_{i,\varepsilon}) \cap W_{j,\varepsilon'}$ is a wall in the rank $r-1$ median algebra $L(W_{i,\varepsilon})$, we deduce by Lemma~\ref{lem:rank_lw} that $L(W_{i,\varepsilon}) \cap L(W_{j,\varepsilon'})$ has rank $r-2$.

By induction, we prove that for each $p \in \llb 1,r\rrb$, for each distinct $i_1,\ldots,i_p \in \llb 1,r\rrb$ and each $\varepsilon_1,\ldots,\varepsilon_p \in \{0,1\}$, the intersection $\bigcap_{1 \leq k \leq p} L(W_{i_k,\varepsilon_k})$ has rank $r-p$.

\bigskip

For each $k \in \llb 1,r \rrb$, let $e_k=(0,\ldots,0,1,0,\ldots,0)$ (the $1$ is in the $k^{th}$ position). Hence for each $k \in \llb 1,r \rrb$, the points $f(0)$ and $f(e_k)$ are contained in a convex rank $1$ closed subspace. In particular, there exists an injective, monotone, biLipschitz path $p_k$ from $f(0)$ to $f(e_k)$, with convex image.

We will show by induction on $k \in \llb 0,r\rrb$ that we can extend $f$ to a biLipschitz median embedding from $[0,1]^k \times \{0,1\}^{r-k}$ into $X$. The case $k=0$ is already true.

\begin{figure}[!h]
\def\svgwidth{6cm}
\center
\input{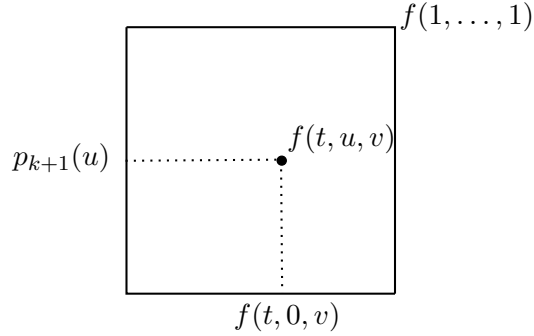}
\caption{Extending f}
\label{fig:Convex_cube}
\end{figure}

Assume we have extended $f$ to a biLipschitz median embedding $f:[0,1]^k \times \{0,1\}^{r-k} \ra X$ for some $k<r$. Define
\beq f:[0,1]^k \times [0,1] \times \{0,1\}^{r-k-1} &\ra& X \\
(t,u,v) \in [0,1]^k \times [0,1] \times \{0,1\}^{r-k-1} & \mapsto & \mu(f(t,0,v),p_{k+1}(u),f(1,\ldots,1)).\eeq
See Figure~\ref{fig:Convex_cube}. Since $p_{k+1}$ and $\mu$ are biLipschitz, we deduce that $f$ is biLipschitz on $[0,1]^k \times [0,1] \times \{0,1\}^{r-k-1}$.

So $f$ is extended to a biLipschitz map $f:[0,1]^r \ra X$. If $t \in [0,1]^r$ and $k \in \llb 1,r\rrb$, notice that in $[0,1]^r$ the median of $(t,0,e_k)$ is equal to $t_ke_k$. So $\mu(f(t),f(0),f(e_k))=f(t_ke_k)=p_k(t_k)$. Since each path $p_1,\ldots,p_r$ is injective, we deduce that $f$ itself is injective. 

By using several times Property~(\ref{eq:abcde}), we prove that $f$ preserves the medians. Hence $f$ is a median embedding, and by Lemma~\ref{lem:convex cube} the image of $f$ is convex.
\ep

Let us recall Bowditch's definition of the separation dimension of a space, which is a good notion of dimension when working with medians on a metric space.

\bdf [Bowditch] If $X$ is a Hausdorff topological space, define the \emph{separation dimension} of $X$ inductively as follows:
\bit
\item If $X=\emptyset$, then the separation dimension of $X$ is $-1$.
\item $X$ has separation dimension at most $n \in \N$ if for any distinct $x,y \in X$, there exist closed subsets $A$, $B$ of $X$ such that $x \not\in B$, $y \not\in A$, $X=A \cup B$ and $A \cap B$ has dimension at most $n-1$.
\eit
\edf

\brk If $X$ is a Hausdorff metric space, then the inductive dimension of $X$ is at most equal to the separation dimension. Conversely, we have the following. \erk 

\blem \cite[Section~III.6]{hurewicz_wallman} If $X$ is a locally compact Hausdorff metric space, then the inductive dimension of $X$ equals the separation dimension. \elem

The following Lemma is immediate.

\blem \label{lem:dimension_injection} Let $X,Y$ be Hausdorff topological spaces, and let $f:X \ra Y$ be a continuous injective map. Then the separation dimension of $X$ is at most equal to the separation dimension of $Y$. \elem

We deduce the following.

\bcor \label{cor:dimension_equals_rank} Let $X$ be a connected metric space, with a Lipschitz locally convex median of rank $r$. Then the separation dimension of $X$ equals $r$. \ecor

\bp According to~\cite[Theorem~2.2]{bowditch_coarse_median}, the separation dimension of $X$ is bounded above by $r$. According to~\ref{pro:embedding_cube}, there exists an embedding of $[0,1]^r$ into $X$, so according to Lemma~\ref{lem:dimension_injection} the separation dimension of $X$ is precisely equal to $r$. \ep

Finally, for affine buildings, we have the following.

\bcor \label{cor:embedding_cube_building} Let $X$ be an affine building of rank $r$. Then any locally convex Lipschitz median on $X$ has rank $r$. \ecor

\bp According to~\cite[Theorem~B]{kramer_local_structure}, $X$ has separation dimension equal to $r$. According to Corollary~\ref{cor:dimension_equals_rank}, any locally convex Lipschitz median on $X$ has rank $r$. \ep

\section{Proof of Theorem~\ref{thm:no median on Xinfty}}

\label{sec:gate}

Consider a thick affine building $X$. Assume that there exists a $k$-Lipschitz, locally convex median $\mu$ on $X$. We will show that $X$ has spherical type $A_1^r$.

\bpro \label{pro:locally_l1} There exists $x \in X$ such that in a tangent cone $(X_\infty,d_\infty,x_\infty,\mu_\infty)$ of $(X,d,x,\mu)$ at $x$, the ultralimit $F_\infty$ of some apartment $F$ containing $x$ is convex and median-isomorphic to $(\R^r,L^1)$ by an affine isomorphism. \epro

\bp According to Corollary~\ref{cor:embedding_cube_building}, the median $\mu$ has rank $r$. According to Proposition~\ref{pro:embedding_cube}, there exists a biLipschitz, median embedding $f$ of $[0,1]^r$ into $X$ with convex image. According to~\cite[Corollary~6.2.3]{kleiner_leeb}, the image of $f$ intersects finitely many apartments of $X$. Consider a non-empty open subset $U$ of $[0,1]^r$ such that $f(U)$ lies in one apartment $F$ of $X$. The map $f|_U : U \ra F$ is biLipschitz, hence it is differentiable almost everywhere: pick a point $t \in U$ where $f$ is differentiable. Since $f$ is biLipschitz, the differential of $f$ at $t$ is invertible. Then in any tangent cone of $(X,d,x,\mu)$ at $x=f(t)$, the ultralimit of $F$ is convex and median-isomorphic to $(\R^r,L^1)$, by an affine isomorphism.\ep

According to Proposition~\ref{pro:locally_l1}, up to considering a tangent cone of $X$ and using Theorem~\ref{thm:ultralimit of affine building}, we can assume that there exists a convex apartment $F$ of $X$ with a median, affine isomorphism with $(\R^r,L^1)$. Since $F$ is convex, closed and locally compact, we can consider $\pi_F : X \ra F$ the gate projection onto $F$.

\blem \label{lem:gate in intersection} For each $x \in X \bs F$, and for each apartment $F'$ of $X$ containing $x$ such that $F \cap F'$ is a half-apartment, we have $\pi_F(x) \in F \cap F'$. \elem

\bp  By contradiction, assume that there exists such an $x \in X \bs F$ and an apartment $F'$ containing $x$ such that $F \cap F'$ is a half-apartment, and such that $\pi_F(x) \not\in F \cap F'$. Fix a Lipschitz embedding $\iota$ of the $(r-1)$-ball $\B^{r-1}$ into $F \cap F'$. Extend $\iota$ to a Lipschitz embedding of the half $r$-ball $\B^{r,+}$ into $F' \bs \overset{\circ}{F}$, where $\B^{r-1}$ is the equatorial sphere of $\B^r$. Extend $\iota$ to a Lipschitz map $\iota : \B^r \ra F \cup F'$ by setting, for $z \in \B^{r,-}$, $\iota(z) = \pi_F(\iota(-z))$. Since $\pi_F(x) \in F \bs F'$ and $\iota$ is Lipschitz, we deduce that $\left(\iota(\B^r) \bs \iota(\partial \B^r) \right) \cap F$ has non-empty interior.

For each $z \in (\partial \B^r)^+ = \SS^{r-1,+}$, we have $\iota(-z)=\pi_F(\iota(z))$. Consider the following map:
\beq \widetilde{\iota'} : \SS^{r-1,+} \times [0,1] & \ra & X \\
(z,t) & \mapsto & \mu(\iota(z),\pi_F(\iota(z)),[(1-t)\iota(z)+t\pi_F(\iota(z))]),\eeq
where $[(1-t)\iota(z)+t\pi_F(\iota(z))]$ is the unique point on the CAT(0) geodesic segment between $\iota(z)$ and $\pi_F(\iota(z))$ at distance $td(\iota(z),\pi_F(\iota(z)))$ from $\iota(z)$ (see Figure~\ref{fig:Sphere}).

\begin{figure}[!h]
\def\svgwidth{8cm}
\center
\input{Sphere.pdf_tex}
\caption{The sphere $\SS^r$ in $X$}
\label{fig:Sphere}
\end{figure}

The map $\widetilde{\iota'}$ is Lipschitz, and it is such that for every $z \in \partial \SS^{r-1,+} = \partial \B^{r-1}$, as $\iota(z) \in F \cap F'$ we have $\pi_F(\iota(z))=\iota(z)$, so for every $t,t' \in [0,1]$, we have $\widetilde{\iota'}(z,t)=\widetilde{\iota'}(z,t')$. Consider the quotient of $\SS^{r-1,+} \times [0,1]$ by the equivalence relation defined for every $z \in \partial \SS^{r-1,+}$ and $t,t' \in [0,1]$ by $(z,t) \sim (z,t')$: it is a topological ball $\B^r$. So $\widetilde{\iota'}$ induces a Lipschitz map $\iota' : \B^r \ra X$ such that $\iota|_{\SS^{r-1}}=\iota'|_{\SS^{r-1}}$. This defines a Lipschitz map $\alpha: \SS^r \ra X$.

\medskip

In $\alpha(\SS^r)$, if we collapse the complement of a small open ball in $F \bs F'$, we obtain a topological sphere $\SS^r$. As a consequence, $H_r(\alpha(\SS^r)) \neq 0$. According to~\cite[Theorem~B]{kramer_local_structure}, $X$ has topological dimension $r$, and since $\alpha(\SS^r)$ is a compact subspace of $X$, we deduce that $H_r(\alpha(\SS^r)) \ra H_r(X)$ is an injection (see for instance~\cite[Theorem~VIII.3']{hurewicz_wallman}). Since $X$ is contractible, this is a contradiction. \ep

We can now conclude the proof of Theorem~\ref{thm:no median on Xinfty}. By contradiction, assume that $X$ has not spherical type $A_1^r$. Since $X$ is thick, there exists a Weyl wall $W$ in $F$, and two singular geodesics $\gamma,\gamma'$ in $X$, each interesecting $W$, such that $\gamma$ and $\gamma'$ intersect in $X \bs F$. Let $x = \gamma \cap \gamma' \in X \bs F$. Since $\gamma$ is singular, the intersection of all apartments $F'$ containg $\gamma$ such that $F' \cap F$ is a half-apartment is equal to $\gamma$. According to Lemma~\ref{lem:gate in intersection}, we deduce that $\pi_F(x) \in \gamma \cap F$. Similarly, $\pi_F(x) \in \gamma' \cap F$. This contradicts the assumption that $\gamma$ and $\gamma'$ intersect in $X \bs F$.

As a consequence, $X$ has spherical type $A_1^r$. This concludes the proof of Theorem~\ref{thm:no median on Xinfty}, as well as Theorem~\ref{thm:no_cm_on_X_maindirection} and Theorem~\ref{thm:no_cm_on_X}.

\section{Proof of the main consequences}

\label{sec:corollaries}

We will now prove the main consequences of Theorem~\ref{thm:no_cm_on_X}, namely Theorem~\ref{thm:not_qi_CCC} and Theorem~\ref{thm:no_cm_on_G}, and also give the proof of Proposition~\ref{pro:rank1}.

\bp[Proof of Theorem~\ref{thm:not_qi_CCC}]
In one direction, assume that $X$ is a symmetric space or affine building of spherical type $A_1^r$. If $X$ is a discrete affine building of spherical type $A_1^r$, if we endow it with the $L^1$ metric, $X$ becomes an actual CAT(0) cube complex. If $X$ is a symmetric space, it is isometric to a product of rank $1$ symmetric spaces. According to~\cite[Theorem~1.8]{haglund_wise}, every word-hyperbolic group is quasi-isometric to CAT(0) cube complex. So each rank $1$ symmetric space is quasi-isometric to a CAT(0) cube complex, hence $X$ itself is quasi-isometric to a CAT(0) cube complex.

\medskip

Conversely, assume that the symmetric space or affine building $X$ is quasi-isometric to a CAT(0) cube complex $(Y,d_p)$, possibly of infinite dimension, endowed with the $L^p$ distance for some $p \in [1,\infty]$. Since $(Y,d_p)$ is quasi-isometric to the metric space $X$ which has finite dimension, we deduce that $(Y,d_p)$ is quasi-isometric to $(Y,d_1)$. Since $(Y,d_1)$ is a metric median space, we deduce that there exists a coarse median on $X$. According to Theorem~\ref{thm:no_cm_on_X}, we deduce that the spherical type of $X$ is $A_1^r$. 
\ep

\bp[Proof of Theorem~\ref{thm:no_cm_on_G}]
Assume that $\Gamma$ is coarse median. Since non-uniform lattices do not have property (RD), $\Gamma$ is cocompact in $G$. So $\Gamma$, endowed with a word metric, is quasi-isometric to $G$, endowed with a left $G$-invariant metric. Let $X$ be the symmetric space of non-compact type of $G$ (if $\K=\R$ or $\C$) or the Bruhat-Tits Euclidean building of $G$ (if $\K$ is non-archimedean). Then $G$ is quasi-isometric to $X$, and so $X$ has a coarse median. According to Theorem~\ref{thm:no_cm_on_X}, $X$ has spherical type $A_1^r$, so $G$ has relative type $A_1^r$. Since $G$ is simple, $r=1$, and $G$ has $\K$-rank $1$.
\ep

We will now consider the rank $1$ case, and give the proof of Proposition~\ref{pro:rank1}.

\bp[Proof of Proposition~\ref{pro:rank1}]
\bit
\item If $\Gamma$ is a uniform lattice in $G$, then $\Gamma$ is hyperbolic hence coarse median.
\item If $G$ is locally isomorphic to $\SO_0(n,1)$ for some $n \geq 2$, $\Gamma$ is hyperbolic relatively to a family $P_1,\dots,P_m$ of parabolic subgroups. Each parabolic subgroup $P_i$ is virtually isomorphic to $\Z^{n-1}$. In particular, each $P_i$ is coarse median, so by~\cite{bowditch_rh_coarse_median} $\Gamma$ itself is coarse median.
\item If $G$ is locally isomorphic to $\SU(2,1)$, $\Gamma$ is hyperbolic relatively to a family $P_1,\dots,P_m$ of parabolic subgroups. Each parabolic subgroup $P_i$ is virtually isomorphic to the $3$-dimensional Heisenberg group $H_3$, which has cubic Dehn function (see~\cite{dehn_heisenberg}). This implies that $\Gamma$ has cubic Dehn function (see~\cite{osin_rh}), so by~\cite{bowditch_embedding_trees} $\Gamma$ is not coarse median.
\eit
\ep

\bibliographystyle{abbrv}
\bibliography{../../../bibli}

\sign

\end{document}